\documentclass[11pt]{amsart}
\usepackage[T1]{fontenc}
\usepackage{lmodern,amsmath,amssymb,mathtools,microtype}
\usepackage[margin=1.1in]{geometry}
\usepackage[hidelinks]{hyperref}
\hypersetup{pdftitle={The Kneser--Poulsen conjecture for d+5 balls and uniform contractions},pdfauthor={Andrii Arman, Andriy Bondarenko, Andriy Prymak, Danylo Radchenko}}
\newtheorem{theorem}{Theorem}[section]
\newtheorem{lemma}[theorem]{Lemma}
\newtheorem{proposition}[theorem]{Proposition}
\newtheorem{corollary}[theorem]{Corollary}
\theoremstyle{definition}

\newtheorem{problem}[theorem]{Open problem}
\theoremstyle{remark}
\newtheorem{remark}[theorem]{Remark}
\numberwithin{equation}{section}
\newcommand{\R}{\mathbb R}
\newcommand{\vol}{\operatorname{vol}}
\newcommand{\ip}[2]{\langle #1,#2\rangle}
\newcommand{\one}{\mathbf 1}

\newcommand{\cU}{\mathcal U}
\newcommand{\cI}{\mathcal I}

\DeclareMathOperator{\diag}{diag}

\title[Kneser--Poulsen conjecture]{The Kneser--Poulsen conjecture for $d+5$ balls and uniform contractions}
\author[A. Arman]{Andrii Arman}
\address{Department of Mathematics, University of Manitoba, Winnipeg, MB, R3T 2N2, Canada}
\email{andrew0arman@gmail.com}
\author[A. Bondarenko]{Andriy Bondarenko}
\address{Department of Mathematical Sciences, Norwegian University of Science and Technology, NO-7491 Trondheim, Norway}
\email{andriybond@gmail.com}
\thanks{The second author was supported in part by Grant 334466 of the Research Council of Norway.}
\author[A. Prymak]{Andriy Prymak}
\address{Department of Mathematics, University of Manitoba, Winnipeg, MB, R3T 2N2, Canada}
\email{prymak@gmail.com}
\thanks{The third author was supported by NSERC of Canada Discovery Grant RGPIN-2026-06488.}
\author[D. Radchenko]{Danylo Radchenko}
\address{Institut des Hautes \'{E}tudes Scientifiques, CNRS, Laboratoire Alexander Grothendieck, 35 route de Chartres, Bures-sur-Yvette 91440, France}
\email{danradchenko@gmail.com}
\thanks{The fourth author was supported by ERC Starting Grant No. 101078782.}
\date{November 29, 2026}

\begin{document}
\begin{abstract}
We prove that a contraction of the centers of at most $d+5$ Euclidean balls in $\R^d$ does not increase the volume of their union and does not decrease the volume of their intersection. The radii may be different. In particular, this proves both inequalities for eight balls in $\R^3$ and nine balls in $\R^4$. We also prove the intersection inequality for uniform contractions of congruent balls, without a restriction on the dimension or the number of balls. More generally, we prove both inequalities for arbitrary radii whenever the endpoint distances can be separated by an intermediate configuration whose boundary spheres have a common point. Both arguments use Csik\'os's first variation formula and differentiation with respect to a common increment of the squared radii. For the cardinality bound, the main step is convexity of the contributions of the walls after integration along a motion in $\R^{d+4}$.
\end{abstract}
\maketitle

\section{Introduction}
Let $P=(p_1,\ldots,p_N)$ and $Q=(q_1,\ldots,q_N)$ be labeled configurations in $\R^d$. We call $Q$ a \emph{contraction} of $P$ if
\begin{equation}\label{eq:contraction}
 |q_i-q_j|\le |p_i-p_j|\qquad(1\le i<j\le N).
\end{equation}
For $r=(r_1,\ldots,r_N)\in[0,\infty)^N$, set
\[
 \cU_d(P,r)=\vol_d\!\left(\bigcup_{i=1}^N B_d(p_i,r_i)\right),\qquad
 \cI_d(P,r)=\vol_d\!\left(\bigcap_{i=1}^N B_d(p_i,r_i)\right),
\]
where $B_d(p,r)$ is the closed Euclidean ball of radius $r$ centered at $p$. The Kneser--Poulsen conjecture asserts that the union and intersection inequalities
\begin{equation}\label{eq:KP}
 \cU_d(Q,r)\le\cU_d(P,r),\qquad
 \cI_d(Q,r)\ge\cI_d(P,r)
\end{equation}
hold whenever \eqref{eq:contraction} holds. The radius $r_i$ remains attached to the $i$th center.

The union inequality for congruent balls was conjectured independently by Poulsen \cite{Poulsen} and Kneser \cite{Kneser}. Csik\'os's first variation formula \cite{Csikos}, in its union and intersection forms, proves both inequalities along continuous contractions. Bezdek and Connelly \cite{BC} combined this formula with motions in a larger Euclidean space and differentiation in the squared radii. They proved \eqref{eq:KP} in the plane for any number of disks, and in $\R^d$ for $N\le d+3$. See \cite{Survey} for further background.

More recently, Wu and Peng \cite{WP} proved both inequalities for $N\le d+4$ in dimensions $d\ge3$. Our first result further increases the cardinality bound by one to $N\le d+5$.
\begin{theorem}\label{thm:cardinality}
Let $d\ge3$, let $N\le d+5$, and let $Q$ be a contraction of $P$ in $\R^d$. Then \eqref{eq:KP} holds for every $r\in[0,\infty)^N$.
\end{theorem}
Thus both inequalities hold for eight balls in three dimensions and nine balls in four dimensions, with arbitrary corresponding radii. In dimensions one and two, both inequalities hold without a cardinality restriction.

Our second result concerns uniform contractions. Following Bezdek and Nasz\'odi \cite{BN}, we say that $Q$ is a \emph{uniform contraction} of $P$ with separating value $a>0$ if
\begin{equation}\label{eq:uniform}
 |q_i-q_j|\le a\le |p_i-p_j|\qquad(i\ne j).
\end{equation}
Here a single value separates every target distance from every source distance; there is no prescribed ratio between corresponding distances.
\begin{theorem}\label{thm:uniform}
Let $d,N\ge1$ and $r>0$. If $Q$ is a uniform contraction of $P$ in $\R^d$, then
\[
 \vol_d\!\left(\bigcap_{i=1}^N B_d(q_i,r)\right)
 \ge
 \vol_d\!\left(\bigcap_{i=1}^N B_d(p_i,r)\right).
\]
\end{theorem}
Bezdek and Nasz\'odi \cite{BN} proved the corresponding comparison for all intrinsic volumes when $N\ge(1+\sqrt2)^d$. Bezdek \cite[Theorem~1]{Bezdek2026} improved this sufficient bound to
\[
 N\ge\left(1+\sqrt{\frac{2d}{d+1}}\right)^d.
\]
Theorem~\ref{thm:uniform} removes the cardinality condition for $d$-dimensional volume. We did not consider the corresponding unrestricted result for all intrinsic volumes.

Theorem~\ref{thm:uniform} follows from the next sufficient condition, which is also of independent interest. It applies to arbitrary radii and has no restriction on the number of balls. The intermediate configuration may lie in a Euclidean space of higher dimension than the two endpoints.
\begin{theorem}\label{thm:intermediate}
Let $P,Q\in(\R^d)^N$ and let $r_i>0$. Suppose that vectors $w_1,\ldots,w_N$ in some Euclidean space satisfy
\begin{equation}\label{eq:intermediate}
 |w_i|=r_i,\qquad
 |q_i-q_j|\le|w_i-w_j|\le|p_i-p_j|\quad(i\ne j).
\end{equation}
Then both inequalities in \eqref{eq:KP} hold.
\end{theorem}
The norm condition means that the spheres $\partial B(w_i,r_i)$ have a common point, chosen as the origin. Thus the theorem applies whenever a contraction can be factored through such a configuration, including when the boundary spheres at either endpoint already have a common point. Condition~\eqref{eq:intermediate} also has an equivalent formulation as a semidefinite feasibility problem; see Remark~\ref{rem:gram}.

For a uniform contraction, the intermediate configuration is a regular simplex whose side length is the separating value $a$. If the source intersection is nonempty, a variance identity shows that this simplex can be placed on a sphere of radius $r$, proving Theorem~\ref{thm:uniform}; the empty-intersection case is immediate. The same construction also proves the union inequality under the radius condition
\[
 r^2\ge\frac{N-1}{2N}a^2,
\]
and hence whenever the source intersection is nonempty; see Corollary~\ref{cor:uniformunion}. For unequal radii, Theorem~\ref{thm:intermediate} remains a sufficient condition for both comparisons. It does not, however, cover arbitrary contractions: the intermediate vectors need not exist, even for congruent balls whose source intersection has nonempty interior. An explicit example is given in Remark~\ref{rem:scope}.

\subsection{Outline of the arguments}
A contraction need not admit a continuous realization in the original space. We use a motion in a larger space for which every squared distance is affine in time. Csik\'os's formula expresses the change of volume as a sum of integrals over the \emph{walls} between truncated Voronoi cells. A common increment $s$ of the squared radii leaves the underlying power diagram unchanged. Differentiation in $s$ then recovers the volume in the original dimension.

For Theorem~\ref{thm:cardinality}, we use four additional coordinates. The desired gain in $d$ dimensions is $1/\pi^2$ times the second derivative of the gain in $d+4$ dimensions. We therefore seek convexity in $s$. For $d+4$ centers this follows from convexity of individual wall volumes at each time. For one more center, we prove convexity after integration in time. In dimensions $d\ge4$, a simultaneous reduction of the two endpoint Gram matrices leads to a weighted integral in two variables. In dimension three, the affine relation among eight paired centers gives a different reduction, according to the signs of its coefficients.

For Theorem~\ref{thm:intermediate}, we compare two motions through the configuration whose boundary spheres have a common point. This common point makes the power diagrams along each motion homothetic after a linear change of coordinates. The intermediate volume cancels before the return to the original dimension. The regular simplex construction then gives Theorem~\ref{thm:uniform}.

\subsection{Related work}
Gorbovickis \cite{Gorbovickis} used higher derivatives in a common squared-radius parameter to prove the union inequality under a restriction on the number of balls meeting each pair intersection in the initial configuration of an expansion. Our cardinality theorem imposes no such intersection condition. Earlier probabilistic volume comparisons under additional assumptions were obtained by Gordon and Meyer \cite{GM}.

Wu and Peng \cite{WP} proved both inequalities for $N\le d+4$ using a motion in $\R^{d+4}$ and convexity of individual wall volumes. Their geometric argument reduces to a three-ball cut by two affine half-spaces. Section~\ref{sec:dplus4} gives an alternative proof using the common radical point and a convex weight on a disk. For $d+5$ centers, our estimates instead involve integration along the motion.

\section{The first variation formula and dimension reduction}\label{sec:framework}
For a configuration $Z=(z_i)$ in $\R^n$ and fixed radii $r_i$, define
\[
 B_i^Z(s)=\{x\in\R^n:|x-z_i|^2\le r_i^2+s\},\qquad s\in\R.
\]
The ball is empty when $r_i^2+s<0$. Write $\cU_n(Z;s)$ and $\cI_n(Z;s)$ for the volumes of the union and intersection of these balls. For the two endpoints in $\R^d$, orient the gains as
\begin{equation}\label{eq:gains}
 g_U(s)=\cU_d(P;s)-\cU_d(Q;s),\qquad
 g_I(s)=\cI_d(Q;s)-\cI_d(P;s).
\end{equation}
Both required inequalities say that the corresponding gain is nonnegative at $s=0$.

These volume functions are continuous in the centers and squared radii. Indeed, the indicators converge away from finitely many spheres, and locally all sets are contained in a fixed bounded set. This also covers balls of radius zero.

\subsection{Power cells and walls}
For centers $z_i\in\R^D$, put
\[
 \phi_i(x)=|x-z_i|^2-r_i^2,
\]
and define the nearest and farthest power cells by
\[
 C_i^- =\{x:\phi_i(x)\le\phi_j(x)\text{ for all }j\},\qquad
 C_i^+ =\{x:\phi_i(x)\ge\phi_j(x)\text{ for all }j\}.
\]
These are also called weighted Voronoi regions. Their truncations by $\phi_i\le s$ partition the union and intersection, respectively, up to sets of measure zero. We use the term \emph{wall}, as in \cite{Csikos,BC}, for a common facet of two such truncated cells. Denote its $(D-1)$-dimensional volume by $W_{ij}^-(t,s)$ or $W_{ij}^+(t,s)$. Only facets between full-dimensional cells are counted.

Let $\ell_{ij}(t)=|z_i(t)-z_j(t)|$. The first variation formulas are
\begin{align}
 \frac{d}{dt}\cU_D(Z(t);s)&=\sum_{i<j}\ell'_{ij}(t)W_{ij}^-(t,s),\label{eq:variationU}\\
 \frac{d}{dt}\cI_D(Z(t);s)&=-\sum_{i<j}\ell'_{ij}(t)W_{ij}^+(t,s).\label{eq:variationI}
\end{align}
We use their integrated form for piecewise analytic motions; see \cite{Csikos,BC}. The opposite signs arise from the opposite orientations of the common facets in the nearest and farthest diagrams.

\subsection{Affine interpolation of squared distances}
The motion
\begin{equation}\label{eq:motion}
 z_i(t)=(\sqrt{1-t}\,p_i,\sqrt t\,q_i),\qquad 0\le t\le1,
\end{equation}
in $\R^d\oplus\R^d$ satisfies
\begin{equation}\label{eq:affinedist}
 \ell_{ij}(t)^2=(1-t)|p_i-p_j|^2+t|q_i-q_j|^2.
\end{equation}
If $N-1\le D$, the same distance data can be realized in $\R^D$. To see this, anchor one label at zero and interpolate the two Gram matrices. On $0<t<1$, the kernel of their positive combination is the intersection of their kernels. The positive square root on the fixed complementary subspace is analytic and gives the required realization. At the endpoints the configurations are congruent to the original configurations embedded in $\R^D$.

Put $\Delta_{ij}=|p_i-p_j|^2-|q_i-q_j|^2$. Integrating \eqref{eq:variationU} or \eqref{eq:variationI} gives the corresponding gain $G(s)$ in dimension $D$:
\begin{equation}\label{eq:integrated}
 G(s)=\frac12\sum_{\Delta_{ij}>0}\Delta_{ij}
        \int_0^1\frac{W_{ij}(t,s)}{\ell_{ij}(t)}\,dt.
\end{equation}
Here and below the argument is applied separately to the two types of walls. These integrals are finite, locally uniformly in $s$. A truncated wall lies in a section of a ball of bounded radius, and, for a shortened pair,
\[
 \ell_{ij}(t)^2\ge(1-t)\Delta_{ij},\qquad
 \int_0^1\ell_{ij}(t)^{-1}\,dt\le\frac2{\sqrt{\Delta_{ij}}}.
\]
One may consequently integrate first over compact subintervals of $(0,1)$ and then pass to the endpoints.

\begin{lemma}[Dimension reduction]\label{lem:slicing}
Embed a configuration in $\R^d$ into $\R^{d+2m}$, where $m\ge1$. For either its union or its intersection volume, denoted by $F_n(s)$ in dimension $n$, one has
\begin{equation}\label{eq:slicing}
 F_{d+2m}(s)=\frac{\pi^m}{(m-1)!}
       \int_{-\infty}^s(s-u)^{m-1}F_d(u)\,du.
\end{equation}
In particular, for either endpoint gain,
\begin{equation}\label{eq:descent}
 G^{(m)}(s)=\pi^m g(s).
\end{equation}
\end{lemma}
\begin{proof}
At a perpendicular coordinate $y\in\R^{2m}$, every squared radius is reduced by $|y|^2$. The fiber is therefore the original union or intersection at parameter $s-|y|^2$. Fubini's theorem and polar coordinates in $\R^{2m}$ give \eqref{eq:slicing}. The function $F_d$ is continuous and vanishes for sufficiently negative arguments, so differentiation gives \eqref{eq:descent}.
\end{proof}

With two extra coordinates, the monotonicity of each wall volume in $s$ suffices: \eqref{eq:integrated} gives $G'\ge0$, and $G'=\pi g$. This recovers the Bezdek--Connelly argument whenever the motion is realizable in $\R^{d+2}$, in particular for $N\le d+3$. For $d=2$, \eqref{eq:motion} already lies in $\R^4$ for every $N$. With four extra coordinates, we instead need convexity of $G$.

\subsection{Changes of coordinates on a wall}
For a nonconstant affine function $f$, the coarea formula gives
\begin{equation}\label{eq:coarea}
 \int_{f=0}\frac{H(x)}{|\nabla f|}\,dS(x)
 =\lim_{\eta\downarrow0}\frac1{2\eta}\int_{|f(x)|<\eta}H(x)\,dx,
\end{equation}
whenever the indicated trace is continuous in the normal coordinate. We use this identity for the ball truncations and polyhedral regions below; changes on their boundaries have zero measure. If $f=\phi_i-\phi_j$, then $|\nabla f|=2\ell_{ij}$. Thus the quotient $W_{ij}/\ell_{ij}$ is twice the integral in \eqref{eq:coarea}. This permits a change of all ambient coordinates using the ordinary determinant, before restriction to the bisector.

\subsection{Approximation of contractions}\label{sec:approximation}
For the cardinality theorem, it is enough to consider generic strict contractions with exactly $d+5$ labels and positive radii. First repeat a labeled ball if necessary; this changes neither volume. Replace $q_i$ by $(1-\eta)q_i$. Every inequality corresponding to distinct source centers becomes strict. Coincident source centers have coincident target centers. A sufficiently small perturbation separating the former makes these inequalities strict while preserving all others. There is now an open set of strict contractions in which both endpoints, and if necessary the radii, can be perturbed into general position. Volume continuity then gives the result for the original configuration. In particular, zero radii are recovered by a limit.

\section{An intermediate configuration with a common boundary point}\label{sec:intermediate}
We prove Theorem~\ref{thm:intermediate} and then deduce its consequences for uniform contractions. The common boundary point makes the power diagrams along each part of the motion homothetic after a linear change of coordinates. The resulting one-variable integral has a nonnegative derivative of the order needed for dimension reduction.

\begin{lemma}\label{lem:scalar}
Let $m\ge1$ be an integer, $\gamma\ge2m-1$, $A>0$, and $b\in\R$. The function
\[
 I(v)=\int_0^1t^\gamma(1-t)^{-m}
       \one_{\{At^2/(1-t)+bt<v\}}\,dt
\]
is smooth on $(0,\infty)$ and satisfies $I^{(m)}(v)\ge0$ there.
\end{lemma}
\begin{proof}
Fix $v>0$. The quadratic polynomial
\[
 Q_v(t)=(1-t)(v-bt)-At^2
\]
has a unique simple root $u\in(0,1)$, is positive on $[0,u)$, and factors as
\[
 Q_v(t)=(u-t)L(t).
\]
Here $L$ is affine and positive on $[0,1]$, since $L(0)=v/u>0$ and $L(1)=A/(1-u)>0$. Thus
\[
 I(v)=\int_0^u t^\gamma(1-t)^{-m}\,dt.
\]
The implicit function theorem gives smoothness in $v$.

Keep $L$ fixed at the chosen value of $v$, and write $u_h$ for the root associated with $v+h$. Then
\[
 u_h=u+h\frac{1-u_h}{L(u_h)}.
\]
The Lagrange inversion formula, applied to the last integral, gives
\begin{equation}\label{eq:lagrange}
 I^{(m)}(v)=(m-1)!\,[x^{m-1}](u+x)^\gamma L(u+x)^{-m}.
\end{equation}
For clarity, the coefficient identity used here is
\[
 [h^n]\int_u^{u+x(h)}f(t)\,dt
 =\frac1n[z^{n-1}]f(u+z)a(u+z)^n,
 \qquad x(h)=h\,a(u+x(h)).
\]
It follows by the formal change of variable $h=x/a(u+x)$ in the residue giving the coefficient; all functions involved are analytic near $u$.

Set $\lambda=uL'/L(u)<1$. After a positive rescaling, the coefficient in \eqref{eq:lagrange} has the sign of
\[
 [z^{m-1}](1+z)^\gamma(1+\lambda z)^{-m}
 =\sum_{j=0}^{m-1}\binom{m+j-1}{j}(1-\lambda)^j
                     \binom{\gamma-m-j}{m-1-j}.
\]
This identity follows by expanding
\[
 (1+z)^\gamma(1+\lambda z)^{-m}
 =\sum_{j\ge0}\binom{m+j-1}{j}(1-\lambda)^j
                         z^j(1+z)^{\gamma-m-j}.
\]
Each summand is nonnegative because $\gamma\ge2m-1$.
\end{proof}

We will use the lemma under positive integrals. If $\Delta_hF(v)=F(v+h)-F(v)$, then
\[
 \Delta_h^m I(v)=\int_{[0,h]^m}I^{(m)}(v+s_1+\cdots+s_m)\,ds_1\cdots ds_m\ge0
 \qquad(v,h>0).
\]
This property is preserved by positive integration and pointwise finite limits. It avoids differentiation under an improper integral.

\begin{proof}[Proof of Theorem~\ref{thm:intermediate}]
We first assume that the vectors $w_i$ are affinely independent. Choose an even integer $k=2m\ge2$ large enough to contain them, and work in $\R^k\oplus\R^d$. For either endpoint $X=(x_i)$, consider
\[
 z_i(t)=(\sqrt{1-t}\,w_i,\sqrt t\,x_i),\qquad 0\le t\le1.
\]
For $X=P$ this parametrizes an expansion; for $X=Q$ it parametrizes a contraction. In each case the gain oriented from $P$ towards $Q$ is a positive sum of the wall integrals in \eqref{eq:integrated}, with the absolute squared-distance change as coefficient.

For the spatial variable $u=(u_W,u_X)$, put
\[
 y=(\sqrt{1-t}\,u_W,\sqrt t\,u_X),\qquad y=t\xi.
\]
Since $|w_i|^2=r_i^2$, the power of the $i$th ball, apart from the common term $|u|^2$, equals
\[
 t(|x_i|^2-r_i^2)-2\ip{(w_i,x_i)}y.
\]
After the substitution $y=t\xi$, all comparisons of powers are independent of $t$. The untruncated wall is therefore represented by a fixed polyhedron $\Pi_{ij}$ in a fixed hyperplane in the $\xi$ coordinates.

Apply \eqref{eq:coarea} to the quotient $W_{ij}/\ell_{ij}$. The first change of variables has determinant factor $(1-t)^{-m}t^{-d/2}$. The dilation contributes $t^{d+k-1}$, including the inverse factor $t^{-1}$ from the affine power difference. Consequently the full factor is
\[
 t^\gamma(1-t)^{-m},\qquad \gamma=2m+d/2-1\ge2m-1.
\]
For a fixed $\xi=(\xi_W,\xi_X)\in\Pi_{ij}$, the truncated power is
\[
 \frac{At^2}{1-t}+bt,
 \quad A=|\xi_W|^2,
 \quad b=|\xi_X|^2-2\ip{(w_i,x_i)}\xi+|x_i|^2-r_i^2.
\]
Thus the time-integrated wall is a positive constant times the integral over $\Pi_{ij}$ of the kernel in Lemma~\ref{lem:scalar}, with $v=s$.

Replace $A$ by $A+\varepsilon$ if necessary. Lemma~\ref{lem:scalar} gives nonnegative $m$th forward differences for $s>0$. Positive integration and then monotone convergence as $\varepsilon\downarrow0$ preserve this property. Finiteness follows from the wall bound after \eqref{eq:integrated}. Summing the gains for the two parts of the motion cancels the intermediate volume at every $s$. The resulting gain $G$ is exactly the gain between the embedded endpoints in dimension $d+2m$. By \eqref{eq:descent}, it is $C^m$ and satisfies $G^{(m)}=\pi^m g$. Nonnegative $m$th forward differences imply $g(s)\ge0$ for $s>0$. Continuity gives the assertion at $s=0$.

To remove affine independence, the following approximation preserves \eqref{eq:intermediate}. For $\varepsilon>0$, put $\delta=\varepsilon^{1/4}$ and replace
\[
 w_i\ \text{by}\ (w_i,\varepsilon e_i),\qquad
 r_i\ \text{by}\ \sqrt{r_i^2+\varepsilon^2},\qquad
 q_i\ \text{by}\ (1-\delta)q_i,
\]
where $e_1,\ldots,e_N$ are new orthonormal vectors. Replace $p_i$ by $(1+\delta)p_i+\delta^2v_i$, where the $v_i\in\R^d$ are fixed and pairwise distinct. For originally distinct source centers, the expansion creates a margin of order $\delta$, which dominates the perturbation of order $\delta^2$ and the change in intermediate distances. If $p_i=p_j$, then $w_i=w_j$ and $q_i=q_j$, while the new source distance is $\delta^2|v_i-v_j|>\sqrt2\,\varepsilon$ for sufficiently small $\varepsilon$. Hence both distance inequalities hold. The new intermediate vectors are affinely independent and have the new prescribed norms. Apply the proved case and let $\varepsilon\downarrow0$.
\end{proof}

\begin{remark}\label{rem:gram}
Condition \eqref{eq:intermediate} is equivalent to the existence of a positive semidefinite matrix $M$ such that
\[
 M_{ii}=r_i^2,\qquad
 \frac{r_i^2+r_j^2-|p_i-p_j|^2}{2}
 \le M_{ij}\le
 \frac{r_i^2+r_j^2-|q_i-q_j|^2}{2}.
\]
Indeed, $M$ is the Gram matrix of the intermediate vectors. The condition can therefore be tested by semidefinite feasibility.
\end{remark}

\begin{remark}[Scope of the intermediate-configuration condition]\label{rem:scope}
Condition~\eqref{eq:intermediate} need not hold for a strict contraction, even when all radii are equal and the source intersection has nonempty interior. Let $e\in\R^d$ be a unit vector and take
\[
 P=(-e/2,0,e/2),\qquad Q=\frac{49}{50}P,\qquad r_1=r_2=r_3=1.
\]
The source intersection contains $B_d(0,1/2)$. If intermediate vectors existed, their unit norms and the bounds $|w_1-w_2|,|w_3-w_2|\le1/2$ would imply
\[
 \ip{w_1}{w_2},\ \ip{w_3}{w_2}\ge\frac78.
\]
For $m=(w_1+w_3)/2$, it follows that $|m|\ge\ip{m}{w_2}\ge7/8$. Hence
\[
 |w_1-w_3|^2=4(1-|m|^2)\le\frac{15}{16}
 <\left(\frac{49}{50}\right)^2=|q_1-q_3|^2,
\]
contrary to \eqref{eq:intermediate}. Thus the sufficient condition does not cover all contractions and does not by itself settle the conjecture in general.
\end{remark}

\subsection{Uniform contractions}
For a uniform contraction, the separating value is the side length of the intermediate regular simplex. The only issue is whether that simplex can be placed on the prescribed common sphere. A point in the source intersection gives exactly the required circumradius bound.
\begin{proof}[Proof of Theorem~\ref{thm:uniform}]
The case $N=1$ is immediate. If the source intersection is empty, the assertion is also immediate. Otherwise choose $x\in\bigcap_iB_d(p_i,r)$ and put $\bar p=N^{-1}\sum_i p_i$. The variance identity gives
\[
 Nr^2\ge\sum_i|p_i-x|^2\ge\sum_i|p_i-\bar p|^2
 =\frac1N\sum_{i<j}|p_i-p_j|^2
 \ge\frac{N-1}{2}a^2.
\]
Thus $\rho_N^2=(N-1)a^2/(2N)\le r^2$. Let $v_i$ be the vertices of a regular $(N-1)$-simplex of side $a$, centered at zero, so $|v_i|=\rho_N$. Choose a unit vector $n$ perpendicular to its span and set
\[
 w_i=v_i+\sqrt{r^2-\rho_N^2}\,n.
\]
Then $|w_i|=r$ and $|w_i-w_j|=a$. Theorem~\ref{thm:intermediate} applies.
\end{proof}

\begin{corollary}\label{cor:uniformunion}
For a uniform contraction with separating value $a$, both volume inequalities hold for congruent balls of radius $r$ whenever
\[
 r^2\ge\frac{N-1}{2N}a^2.
\]
In particular, both hold if the source intersection is nonempty.
\end{corollary}
\begin{proof}
The same regular simplex satisfies Theorem~\ref{thm:intermediate}.
\end{proof}
The empty-intersection case used in Theorem~\ref{thm:uniform} gives no conclusion for the union. Thus Corollary~\ref{cor:uniformunion} retains its radius hypothesis.

\section{The case of \texorpdfstring{$d+4$}{d+4} balls}\label{sec:dplus4}
We first isolate a simple convexity observation. It explains the role of the free coordinate in the $d+4$ case and will also be used for eight balls in $\R^3$.
\begin{lemma}\label{lem:disk}
If $\psi:\R\to[0,\infty)$ is finite and convex, then
\[
 H(v)=\int_{x^2+y^2\le v}\psi(x)\,dx\,dy\quad(v>0),
 \qquad H(v)=0\quad(v\le0),
\]
is convex on $\R$.
\end{lemma}
\begin{proof}
Polar coordinates give
\[
 H'(v)=\frac12\int_0^{2\pi}\psi(\sqrt v\cos\theta)\,d\theta\qquad(v>0).
\]
For each $c\ge0$, the even convex function $u\mapsto\psi(uc)+\psi(-uc)$ is nondecreasing for $u\ge0$. Pairing opposite angles proves that $H'$ is nondecreasing. Its right limit at zero is $\pi\psi(0)\ge0$, which proves convexity across zero as well.
\end{proof}

\begin{proposition}\label{prop:dplus4}
Both inequalities in \eqref{eq:KP} hold for $N\le d+4$ balls in $\R^d$, $d\ge3$.
\end{proposition}
\begin{proof}
By the approximation in Section~\ref{sec:approximation}, with $d+4$ labels, we may assume that the interior configurations of \eqref{eq:affinedist} have affine dimension $d+3$. Realize the motion in $\R^{d+4}$. At a fixed time, let $A$ be the affine span of the centers. There is a unique point $o\in A$ at which all powers are equal: subtracting the powers gives $d+3$ independent affine equations. Thus the nearest and farthest power cells in $A$ are polyhedral cones with vertex $o$.

Fix a wall. Its part in $A$ has dimension $d+2$. Polar coordinates about $o$ write its measure as $\rho^{d+1}\,d\rho\,d\omega$, with $\rho\ge0$ and $\omega$ in a fixed set of directions. There is also one free coordinate perpendicular to $A$. For a fixed direction $\omega$, the ray and this perpendicular coordinate lie in a plane, and the ball truncation is a disk with squared radius $s$ plus a constant. After translating the coordinate along the ray, its weight is $(c+x)_+^{d+1}$, a nonnegative convex function. Lemma~\ref{lem:disk} proves convexity in $s$ of the wall volume at this fixed time. Integrating over $\omega$ and then in time preserves convexity. Equations~\eqref{eq:integrated} and~\eqref{eq:descent}, with $m=2$, give the result.
\end{proof}

\begin{remark}[Comparison with Wu--Peng]\label{rem:WP}
Wu and Peng prove convexity of $s\mapsto\vol_m(K\cap B_m(c,\sqrt s))$ when $m\ge3$ and $K$ is cut out by at most $m-1$ affine half-spaces, taking the value zero for $s<0$. They use an angular overlap formula and slicing. For $N=d+5$, a wall has $d+3$ inequalities in dimension $d+3$, beyond this half-space count. We prove convexity of the time-integrated contributions instead. Both arguments use finite differences to pass to limits before recovering the original-dimensional volumes.
\end{remark}

\section{At most \texorpdfstring{$d+5$}{d+5} balls in dimensions \texorpdfstring{$d\ge4$}{d at least 4}}\label{sec:highdim}
For one additional center, the interior configurations can fill $\R^{d+4}$, so the perpendicular coordinate used in Proposition~\ref{prop:dplus4} is no longer available. The time variable takes its place. The following planar integral is the analytic ingredient; its proof is given in Appendix~\ref{sec:planarproof}.

\begin{lemma}\label{lem:planar}
Let $Y_P,Y_Q$ be vector-valued affine functions of $(t,\rho)$, and let $Y_1,\ldots,Y_k$ be scalar-valued affine functions of $(t,\rho)$. Let $l_1,\ldots,l_k$ be affine functions of $t$, positive on $[0,1]$. Put
\begin{align}
 B(t)&=t(1-t),\notag\\
 F(t,\rho)&=\frac{|Y_P(t,\rho)|^2}{1-t}
       +\frac{|Y_Q(t,\rho)|^2}{t}
       +\sum_{j=1}^k\frac{Y_j(t,\rho)^2}{l_j(t)},\label{eq:planarF}\\
 h(t,\rho)&=\rho_+^m\prod_{j=1}^k l_j(t)^{-1/2},\notag
\end{align}
where $m\ge4+k/2$ is an integer. Then
\begin{equation}\label{eq:planarH}
 H(v)=\int_0^1\int_\R h(t,\rho)B(t)^{-2}
                      \one_{\{F(t,\rho)<v\}}\,d\rho\,dt
\end{equation}
is convex as an extended-valued function. In particular, it is convex on each interval on which it is finite.
\end{lemma}
The singular weight $B^{-2}$ is canceled by two differentiations of the moving level boundary of $B(v-F)$. This leaves a weighted area in the $(t,\rho)$-plane. Convexity follows by comparing the two ends of each chord through the maximum of a regularized level function.

\begin{proof}[Proof of Theorem~\ref{thm:cardinality} for $d\ge4$]
By Section~\ref{sec:approximation}, assume that $N=d+5$, both endpoint configurations have affine rank $d$, and their paired configuration has affine rank $d+4$. Anchor one label at zero and let $M_P,M_Q$ be the $(d+4)\times(d+4)$ Gram matrices of the remaining vectors. Then $M_P+M_Q$ is positive definite.

A congruence making $M_P+M_Q$ the identity, followed by an orthogonal diagonalization, gives an invertible matrix $E$ such that
\[
 M_P=E^{\mathsf T}\diag(\theta_1,\ldots,\theta_{d+4})E,\qquad
 M_Q=E^{\mathsf T}\diag(1-\theta_1,\ldots,1-\theta_{d+4})E.
\]
Both forms have nullity four and their kernels are disjoint. Hence four of the $\theta_j$ equal $1$, four equal $0$, and the remaining $k=d-4$ lie in $(0,1)$. Let $e_i$ be the columns of $E$, together with the anchored zero. Define $S_t$ to be diagonal with entries
\[
 \underbrace{\sqrt{1-t},\ldots,\sqrt{1-t}}_{4},\quad
 \underbrace{\sqrt t,\ldots,\sqrt t}_{4},\quad
 \sqrt{l_1(t)},\ldots,\sqrt{l_k(t)},
\]
where $l_j(t)=(1-t)\theta_j+t(1-\theta_j)>0$. The centers $S_te_i\in\R^{d+4}$ realize \eqref{eq:affinedist}.

Fix a shortened pair $i,j$ and translate its first site to zero, replacing $e_l$ by $e_l-e_i$. For a spatial point $x$, put $y=S_tx$. The power differences are
\begin{equation}\label{eq:slacksHD}
 f_l(y,t):=\phi_i(x)-\phi_l(x)
       =2\ip{e_l}y-|S_te_l|^2+r_l^2-r_i^2.
\end{equation}
Their linear parts in $y$ are fixed and their offsets are affine in $t$. On the bisector $f_j=0$, the remaining $d+3$ differences are independent affine coordinates. Choose their common sign according to the type of power cell and denote the resulting nonnegative coordinates by $h_l$. Thus the wall corresponds to $h_l\ge0$.

Write these coordinates in the form
\[
 h_l=\rho\alpha_l,\qquad \rho\ge0,\quad \alpha_l\ge0,\quad\sum_l\alpha_l=1.
\]
The Jacobian in the $(d+3)$-dimensional orthant is a positive constant times $\rho^{d+2}$. For each fixed $\alpha$, solving \eqref{eq:slacksHD} makes $y$ affine in $(t,\rho)$. Splitting $y$ into the four source coordinates, four target coordinates, and $k$ mixed coordinates gives
\[
 |x|^2=|S_t^{-1}y|^2=F(t,\rho)
\]
with $F$ as in \eqref{eq:planarF}. The ball truncation is $F\le r_i^2+s$.

By \eqref{eq:coarea}, the factor $\ell_{ij}^{-1}$ is absorbed before changing the ambient coordinates. The remaining determinant factor is
\[
 (\det S_t)^{-1}=B(t)^{-2}\prod_{j=1}^k l_j(t)^{-1/2}.
\]
The linear parts of the power coordinates are independent of $t$, so their Jacobian is a positive constant independent of $t$ and $s$. It follows that
\[
 s\longmapsto\int_0^1\frac{W_{ij}(t,s)}{\ell_{ij}(t)}\,dt
\]
is a positive integral, over $\alpha$, of functions \eqref{eq:planarH} evaluated at $v=r_i^2+s$, with
\[
 m=d+2,\qquad k=d-4,\qquad m\ge4+k/2.
\]
Strict and non-strict ball truncations have the same measure here. Lemma~\ref{lem:planar}, positive integration, and the finiteness following \eqref{eq:integrated} prove convexity. Thus $G$ in \eqref{eq:integrated} is convex. Equation~\eqref{eq:descent} with $m=2$ yields $\pi^2g=G''\ge0$. Approximation completes the proof.
\end{proof}

\section{Eight balls in three dimensions}\label{sec:three}
For eight centers in $\R^3$, the paired centers $(p_i,q_i)$ lie in $\R^6$ and have one affine relation in general position. A wall in $\R^7$ has one free perpendicular coordinate. Its other coordinates can be described by six nonnegative numbers subject to one affine relation. If one sign group has at least four members, a simplex of dimension at least three supplies the power needed for the following time-integrated estimate. The remaining $3+3$ case admits a fixed-time argument using Lemma~\ref{lem:disk}.

\begin{lemma}\label{lem:interval}
Let $a,b>0$, let $A,C,k\in\R$, and let $m\ge3$ be an integer. Then
\begin{equation}\label{eq:intervalJ}
 J(v)=\int_0^1\frac{(A+Ct)_+^m}{[t(1-t)]^{3/2}}
       \left(v-\frac a{1-t}-\frac bt+kt\right)_+^{1/2}\,dt
\end{equation}
is convex on $\R$.
\end{lemma}
The square-root factor comes from the free coordinate. Although differentiating it twice at fixed $t$ is not useful, transporting the two endpoints of the contributing time interval gives a positive expression.
\begin{proof}
Put $w(t)=(A+Ct)_+^m$ and $\varphi(t)=a/(1-t)+b/t-kt$. The function $\varphi$ is strictly convex and diverges at both endpoints. For $v>\min\varphi$, its strict sublevel set is an interval $(l,u)\subset(0,1)$. The cubic polynomial
\[
 Q_v(t)=t(1-t)(v-\varphi(t))
\]
factors as $Q_v(t)=T(t)D(t)$, where
\[
 T(t)=(t-l)(u-t),\qquad L=u-l,
\]
and $D$ is affine and positive on $[l,u]$. Positivity at the endpoints follows from the nonzero slopes of $Q_v$. One differentiation gives
\begin{equation}\label{eq:Jprime}
 2J'(v)=\int_l^u w(t)A_v(t)\,dt,
 \qquad A_v(t)=\frac1{t(1-t)\sqrt{Q_v(t)}}.
\end{equation}
To differentiate again, use the velocity
\[
 \eta(v,t)=-\frac{t(1-t)T'(t)}{L^2D(t)}.
\]
Implicit differentiation of $Q_v(l)=Q_v(u)=0$ gives
\[
 l'=-\frac{l(1-l)}{LD(l)},\qquad
 u'=\frac{u(1-u)}{LD(u)}.
\]
These are exactly $\eta(v,l)$ and $\eta(v,u)$. Its local flow therefore maps the old interval to the new one. After this substitution the endpoint square-root singularities are fixed, with a differentiable amplitude. Differentiating \eqref{eq:Jprime} in these coordinates is justified by domination by a constant times $[(t-l)(u-t)]^{-1/2}$.

Set $f=wD^{-3/2}$, with $D$ fixed at the current value of $v$. The combined derivative is integrable, since
\[
 w\partial_vA_v=-\frac{f}{2T^{3/2}},\qquad
 \eta wA_v=-\frac{T'f}{L^2\sqrt T},\qquad
 \left(\frac{T'}{\sqrt T}\right)'=-\frac{L^2}{2T^{3/2}}.
\]
The singular terms cancel, and integration by parts yields
\begin{equation}\label{eq:Jsecond}
 J''(v)=-\frac1{2L^2}\int_l^u\frac{T'}{\sqrt T}f'\,dt
       =\frac1{L^2}\int_l^u\sqrt T\,f''\,dt.
\end{equation}
There are no boundary terms because $\sqrt T$ vanishes at $l,u$.

The function $(x,z)\mapsto x_+^m z^{-3/2}$ is convex for $z>0$. For $x>0$ the diagonal entries of its Hessian are nonnegative and its determinant is
\[
 \frac{3m}{2}\left(m-\frac52\right)x^{2m-2}z^{-5}\ge0.
\]
For $x\le0$ the function is zero, and its first two derivatives match at $x=0$. Restricting to the affine line $(x,z)=(A+Ct,D(t))$ proves $f''\ge0$. Hence \eqref{eq:Jsecond} is nonnegative. Finally, $J$ is continuous and nondecreasing and vanishes on $(-\infty,\min\varphi]$, so convexity extends across the initial level.
\end{proof}

\begin{proof}[Proof of Theorem~\ref{thm:cardinality} for $d=3$]
By Section~\ref{sec:approximation}, assume there are eight labels and that every seven paired points
\[
 e_l=(p_l,q_l)\in\R^6\qquad(1\le l\le8)
\]
are affinely independent. Their unique affine relation is
\begin{equation}\label{eq:affinerelation}
 \sum_{l=1}^8\lambda_le_l=0,\qquad \sum_{l=1}^8\lambda_l=0,
\end{equation}
where every $\lambda_l\ne0$. Use the motion
\[
 z_l(t)=(\sqrt{1-t}\,p_l,\sqrt t\,q_l,0)\in\R^7.
\]
For each shortened pair $i,j$, it suffices to prove convexity of
\begin{equation}\label{eq:timewall3}
 s\longmapsto\int_0^1\frac{W_{ij}(t,s)}{\ell_{ij}(t)}\,dt.
\end{equation}

Translate each endpoint so that $p_i=q_i=0$. In the first six spatial coordinates put
\[
 y=(y_P,y_Q)=(\sqrt{1-t}\,x_P,\sqrt t\,x_Q),
\]
and define
\[
 f_l(y,t)=2\ip{e_l}y-(1-t)|p_l|^2-t|q_l|^2+r_l^2-r_i^2.
\]
The bisector is $f_j=0$. The other wall inequalities are $\epsilon f_l\ge0$, where $\epsilon=1$ for farthest cells and $\epsilon=-1$ for nearest cells. Formula~\eqref{eq:coarea} and integration in the seventh coordinate give
\begin{align}\label{eq:wall3}
 \frac{W_{ij}(t,s)}{\ell_{ij}(t)}
 =\frac4{[t(1-t)]^{3/2}}\int_{f_j=0}
 &\prod_{l\ne i,j}\one_{\{\epsilon f_l\ge0\}}\notag\\[-2mm]
 &\times\left(r_i^2+s-\frac{|y_P|^2}{1-t}-\frac{|y_Q|^2}{t}\right)_+^{1/2}
       \frac{dS_y}{|\nabla_y f_j|}.
\end{align}
The factor $4$ is the product of $2$ from the bisector normal and $2$ from the free interval.

On the five-dimensional hyperplane $f_j=0$, introduce the six nonnegative coordinates
\[
 h_l=|\lambda_l|\epsilon f_l,\qquad l\ne i,j.
\]
Their unique relation, obtained from \eqref{eq:affinerelation}, is
\begin{equation}\label{eq:slackrelation}
 \sum_{l\in L}h_l-\sum_{l\in M}h_l=C(t),
\end{equation}
where $C$ is affine and $L,M$ are the positive and negative sign groups of the remaining $\lambda_l$. Reverse all signs if necessary so that $|L|\ge|M|$. Any five of the six slacks give affine coordinates on the bisector. The corresponding Jacobian, including $|\nabla_y f_j|^{-1}$, is a positive constant independent of $t,s$, because the linear parts of all $f_l$ are fixed and every seven paired points are affinely independent.

Suppose first that $p=|L|\ge4$. Fix the individual slacks $h_M\ge0$ and put
\[
 H(t)=C(t)+\sum_{l\in M}h_l.
\]
For $H(t)>0$, write $h_L=H(t)\alpha$, where $\alpha_l\ge0$ and $\sum_{l\in L}\alpha_l=1$. The simplex has dimension $p-1$, so its Jacobian is a positive constant times $H(t)^{p-1}$. When $H(t)\le0$, its contribution is zero. The factor is therefore $H(t)_+^{p-1}$.

For fixed $h_M,\alpha$, the solution $y(t)$ of the slack equations is affine in $t$, including its formal continuation where $H(t)<0$. Write $y(t)=y(0)+tV$. Polynomial division gives
\begin{equation}\label{eq:quadraticdivision}
 \frac{|y_P(t)|^2}{1-t}+\frac{|y_Q(t)|^2}{t}
 =\frac{|y_P(1)|^2}{1-t}+\frac{|y_Q(0)|^2}{t}+c-kt,
 \qquad k=|V_P|^2-|V_Q|^2,
\end{equation}
for a constant $c$. Substitution in \eqref{eq:wall3} expresses \eqref{eq:timewall3} as a positive integral of the functions in Lemma~\ref{lem:interval}, with
\[
 m=p-1\ge3,\qquad a=|y_P(1)|^2,\qquad b=|y_Q(0)|^2,
 \qquad v=r_i^2+s-c.
\]
If $a$ or $b$ vanishes, replace them by $a+\eta,b+\eta$ and let $\eta\downarrow0$. Monotone convergence and positive integration preserve convexity. Tonelli's theorem applies to all these nonnegative integrals.

It remains to consider $|L|=|M|=3$. Fix $t$, and interchange the groups if necessary so that $C=C(t)\ge0$. Parameterize \eqref{eq:slackrelation} by
\[
 h_L=(C+\rho)\alpha,\qquad h_M=\rho\beta,\qquad \rho\ge0,
\]
where $\alpha,\beta$ range over two-dimensional probability simplices. The measure is a positive constant times
\[
 (C+\rho)^2\rho^2\,d\rho\,d\alpha\,d\beta.
\]
For fixed $\alpha,\beta$, the point in the first six physical coordinates runs along a straight ray. Together with the free seventh coordinate, the line containing it gives a plane in which the truncating ball is a disk. Extend the ray weight by zero to $\rho<0$:
\[
 \psi(\rho)=C^2\rho_+^2+2C\rho_+^3+\rho_+^4.
\]
This is nonnegative and convex. Passage from $\rho$ to Euclidean distance on the line is an affine change, with a positive constant Jacobian. Lemma~\ref{lem:disk} proves convexity in $s$ of this contribution. Positive integration proves convexity of the wall at the fixed time, and hence of \eqref{eq:timewall3}.

Both cases make the gain $G$ in dimension seven convex. Since $G''=\pi^2g$, both three-dimensional comparisons follow. The finiteness needed above follows from \eqref{eq:integrated}. Finally pass to the approximation of Section~\ref{sec:approximation}.
\end{proof}

\section{Further questions}\label{sec:questions}
Theorem~\ref{thm:cardinality} settles the volume inequalities for eight balls in $\R^3$. The next cardinality is therefore nine. We single out the intersection inequality as the next target.
\begin{problem}[Nine balls in three dimensions]\label{prob:nine}
Let $P=(p_1,\ldots,p_9)$ and $Q=(q_1,\ldots,q_9)$ be configurations in $\R^3$ satisfying $|q_i-q_j|\le|p_i-p_j|$ for all $i,j$. Is it true that, for every $r_1,\ldots,r_9>0$,
\begin{equation}\label{eq:nineintersection}
 \vol_3\!\left(\bigcap_{i=1}^9B_3(q_i,r_i)\right)
 \ge
 \vol_3\!\left(\bigcap_{i=1}^9B_3(p_i,r_i)\right)?
\end{equation}
Does the corresponding union inequality also hold?
\end{problem}
Theorem~\ref{thm:uniform} answers \eqref{eq:nineintersection} when the radii are equal and the contraction is uniform. No uniformity or common-boundary hypothesis is imposed in Open problem~\ref{prob:nine}.

The proof for eight balls indicates a specific difficulty in extending the present argument. Nine paired centers in $\R^6$ have two independent affine relations in general position. On a wall, the seven remaining nonnegative slack variables satisfy two affine relations, rather than the single relation \eqref{eq:slackrelation}. The simplex decomposition used in Section~\ref{sec:three} no longer gives the one-variable integral in Lemma~\ref{lem:interval}. A sufficient next step would be convexity of the corresponding time-integrated wall contributions in $\R^7$; convexity of their positively weighted sum would already suffice. These are sufficient approaches, not additional hypotheses in the open problem.

A separate question is whether the intersection conclusion of Theorem~\ref{thm:uniform} remains valid for unequal radii. Theorem~\ref{thm:intermediate} gives a sufficient condition in this direction, but uniform separation of the center distances alone does not supply the intermediate vectors with prescribed norms.

\appendix
\section{Proof of the planar convexity lemma}\label{sec:planarproof}
We prove Lemma~\ref{lem:planar}. Throughout this appendix, $B=t(1-t)$ and $F,h$ are as in its statement.

\subsection{Regularization and the level parameter}
For $\varepsilon>0$, set
\[
 F_\varepsilon=F+\frac\varepsilon B+\varepsilon\rho^2,
 \qquad
 H_\varepsilon(v)=\int_0^1\int_\R hB^{-2}\one_{\{F_\varepsilon<v\}}\,d\rho\,dt.
\]
The square of an affine function divided by a positive affine function is convex: along any line,
\[
 \left(\frac{u^2}{l}\right)''=\frac{2(u'l-ul')^2}{l^3}\ge0.
\]
Hence $F_\varepsilon$ has positive definite Hessian. It is coercive on the strip $(0,1)\times\R$, has a unique minimum $\mu$, and, for $v>\mu$, the sublevel set has compact closure in the strip and smooth boundary.

Fix $v>\mu$ and define
\[
 Q(t,\rho)=B(t)\bigl(v-F_\varepsilon(t,\rho)\bigr),\qquad
 A(c)=\int_0^1\int_\R h(t,\rho)\one_{\{Q(t,\rho)+c>0\}}\,d\rho\,dt.
\]
We first justify the cancellation
\begin{equation}\label{eq:cancelB}
 H_\varepsilon''(v)=A''(0).
\end{equation}
Near the boundary use local coordinates $(r,u)$ with $r=Q$, and a partition of unity. Let $a(r,u)$ be the density of $hB^{-2}\,dt\,d\rho$ in one chart. When $v$ changes to $v+s$, the lower endpoint of the $r$ integral is the root $r_s$ of
\[
 r_s+sB(r_s,u)=0.
\]
At $s=0$, one has $r_0'=-B$ and $r_0''=2B\partial_rB$. The second derivative of the integral is therefore
\[
 -\partial_ra\,(r_0')^2-a r_0''
 =-\partial_r(aB^2)\big|_{r=0}.
\]
For $A(c)$ the density is $aB^2$ and the lower endpoint is $-c$, giving the same expression. Terms supported away from the boundary are locally constant. This proves \eqref{eq:cancelB}.

All boundary differentiations take place in a common compact subset of the strip. Indeed, for $|c|<\varepsilon/2$,
\[
 Q+c\le Bv-\varepsilon-\varepsilon B\rho^2+c.
\]
Positivity forces $Bv>\varepsilon/2$, and then bounds $|\rho|$. This also justifies the use of finitely many boundary charts.

\subsection{The restriction to a line}
On $\{Q>0\}$ the function
\[
 \log Q=\log B+\log(v-F_\varepsilon)
\]
has negative definite Hessian. Thus $Q$ has a unique maximum there. Take this point as origin in the plane and write its value as $z>0$. On a line through it, let $x$ be signed Euclidean distance. We claim that
\begin{equation}\label{eq:gform}
 Q(x)=z-x^2g(x),\qquad
 g(x)=a+bx-cx^2-\sum_j\frac{C_j}{l_j(t(x))},
 \qquad c,C_j\ge0.
\end{equation}
Constants and summands in this representation may depend on the line; terms with a constant denominator can be incorporated into the polynomial part.

To verify the claim, expand
\[
 Q=Bv-t|Y_P|^2-(1-t)|Y_Q|^2
       -\sum_j\frac{BY_j^2}{l_j}-\varepsilon-\varepsilon B\rho^2.
\]
The first three terms have degree at most three on a line. Terms with a constant $l_j$, and the term $-\varepsilon B\rho^2$, have a nonnegative coefficient of $x^4$. For a nonconstant $l_j(t(x))$, polynomial division of $-BY_j^2$ by $l_j$ leaves a nonnegative remainder: at the root $t_*$ of $l_j$ it equals
\[
 -t_*(1-t_*)Y_j(x_*)^2\ge0,
\]
since $t_*\notin[0,1]$. Finally, for a positive affine function $l$,
\[
 \frac1{l(x)}=\frac1{l(0)}-\frac{l'x}{l(0)^2}
                  +\frac{(l')^2x^2}{l(0)^2l(x)}.
\]
Subtract the constant and linear parts of $Q$ and divide by $-x^2$. The linear part vanishes at the chosen maximum, proving \eqref{eq:gform}.

On the closed chord $Q\ge0$, one has $g>0$. At the origin this follows from the negative definite Hessian of $Q$; elsewhere it follows from the strict maximality of $z$. Put
\[
 V=2g+xg',\qquad Q'=-xV.
\]
Strict concavity of $\log Q$ gives $V>0$ in the interior. At each endpoint, the outward radial derivative of $F_\varepsilon$ is positive, by convexity and the strict inequality $F_\varepsilon(0)<v$; hence $V>0$ there as well. Formula~\eqref{eq:gform} implies
\begin{equation}\label{eq:ginequalities}
 g''\le0,\qquad V''=4g''+xg'''\le g'',\qquad 0<\frac Vg<3.
\end{equation}
For the middle inequality, the contribution of $-C/l$ to $V''$ is its contribution to $g''$ multiplied by
\[
 4-\frac{3xl'}l=1+\frac{3l(0)}l\ge1;
\]
the polynomial part has the same sign. The last inequality follows from $V>0$ and the concavity estimate $g-xg'\ge g(0)>0$.

\subsection{The second derivative along a chord}
Use polar coordinates written as unoriented lines through the origin, with area element $|x|\,dx\,d\theta$, $0\le\theta<\pi$. For $c$ sufficiently close to zero, the endpoints $x_-(c)<0<x_+(c)$ of $\{Q+c>0\}$ on each line satisfy
\[
 x_\pm(c)^2g(x_\pm(c))=z+c.
\]
The boundary is compact and regular, and its radial derivatives do not vanish. Consequently these roots and their first two derivatives depend smoothly on $c$, uniformly in the direction. No other components appear for small $|c|$: the sign is bounded away from zero off a neighborhood of the original boundary in a common compact set.

Write $h(x)$ for the restriction of the weight to the line. Since $x_\pm'=1/(x_\pm V(x_\pm))$,
\[
 \frac d{dc}\int_{x_-(c)}^{x_+(c)}|x|h(x)\,dx
 =\frac{h(x_+(c))}{V(x_+(c))}+\frac{h(x_-(c))}{V(x_-(c))}.
\]
Differentiating once more at $c=0$ gives
\begin{equation}\label{eq:chordsecond}
 \left.\frac{d^2}{dc^2}\right|_{c=0}
       \int_{x_-(c)}^{x_+(c)}|x|h(x)\,dx
 =\frac{L(x_+(0))-L(x_-(0))}{\sqrt z},
 \qquad
 L=\frac{\sqrt g\,(h'V-hV')}{V^3}.
\end{equation}
It remains to prove $L'\ge0$ on the chord.

\subsection{The weight inequality}
The function $h^{1/4}$ is convex on the strip. Indeed, for $p,y_j>0$, the monomial
\[
 f(p,y)=p^A\prod_jy_j^{-b_j}
\]
is convex if $b_j\ge0$ and $A\ge1+\sum_jb_j$. Along an affine line, put $X=p'/p$ and $Z_j=y_j'/y_j$. Then
\[
 \frac{f''}{f}=\left(AX-\sum_jb_jZ_j\right)^2
                   +\sum_jb_jZ_j^2-AX^2\ge0.
\]
The inequality follows from Cauchy--Schwarz: the first two terms are at least $A^2X^2/(1+\sum_jb_j)$. Apply this with $A=m/4$, $b_j=1/8$, $p=\rho$, and $y_j=l_j(t)$. The hypothesis $m\ge4+k/2$ gives the required inequality for $A$. Extension by zero to $\rho\le0$ preserves convexity. Hence, wherever $h>0$,
\begin{equation}\label{eq:weight}
 hh''\ge\frac34(h')^2.
\end{equation}

Set
\[
 G=\frac{g'}g,\qquad D=\frac{V'}V,\qquad T=-\frac{V''}V,\qquad R=\frac Vg.
\]
Direct differentiation of $L$ in \eqref{eq:chordsecond} gives
\[
 L'=\frac{\sqrt g}{V^2}
 \left[h''+(G/2-3D)h'+(3D^2-GD/2+T)h\right].
\]
By \eqref{eq:weight}, the bracket is nonnegative if
\begin{equation}\label{eq:discriminant}
 3(3D^2-GD/2+T)-(G/2-3D)^2
 =3T+\frac32GD-\frac14G^2\ge0.
\end{equation}
To verify this, put $U=-g''/g\ge0$ and $\beta=xG<1$. By \eqref{eq:ginequalities},
\[
 R=2+\beta>0,\qquad D=\frac{3G-xU}{R},\qquad T\ge\frac UR.
\]
Therefore
\[
 3T+\frac32GD-\frac14G^2
 \ge\frac{(3-3\beta/2)U+(4-\beta/4)G^2}{R}\ge0.
\]
This proves $L'\ge0$ where $h>0$. At $\rho=0$ the conclusion follows by continuity, since $m\ge4$ and $h,h',h''$ vanish there. On $\rho<0$ they vanish identically.

Integrating \eqref{eq:chordsecond} over the directions gives $A''(0)\ge0$. By \eqref{eq:cancelB}, $H_\varepsilon''(v)\ge0$ for $v>\mu$. The function $H_\varepsilon$ is continuous and nondecreasing and vanishes for $v\le\mu$, so it is convex on $\R$. Finally, $F_\varepsilon\downarrow F$ as $\varepsilon\downarrow0$. Monotone convergence gives $H_\varepsilon(v)\uparrow H(v)$, and the convexity inequality passes to this extended-valued limit. This proves Lemma~\ref{lem:planar}.

\section*{Acknowledgements}
The authors used OpenAI ChatGPT to develop the mathematical arguments and edit the manuscript. The authors take full responsibility for the mathematical content of the work.


\begin{thebibliography}{99}
\bibitem{Bezdek2026}
K.~Bezdek, \emph{Notes on uniform contractions of congruent balls}, preprint, 2026, \href{https://arxiv.org/abs/2602.15352v2}{arXiv:2602.15352v2} (19 October 2026).

\bibitem{BC}
K.~Bezdek and R.~Connelly, \emph{Pushing disks apart---the Kneser--Poulsen conjecture in the plane}, J. Reine Angew. Math. \textbf{553} (2002), 221--236, \href{https://doi.org/10.1515/crll.2002.101}{doi:10.1515/crll.2002.101}.

\bibitem{Survey}
K.~Bezdek, Z.~L\'angi, and M.~Nasz\'odi, \emph{Selected topics from the theory of intersections of balls}, Discrete Appl. Math. \textbf{382} (2026), 60--82, \href{https://doi.org/10.1016/j.dam.2025.11.040}{doi:10.1016/j.dam.2025.11.040}.

\bibitem{BN}
K.~Bezdek and M.~Nasz\'odi, \emph{The Kneser--Poulsen conjecture for special contractions}, Discrete Comput. Geom. \textbf{60} (2018), 967--980, \href{https://arxiv.org/abs/1701.05074}{arXiv:1701.05074}.

\bibitem{Csikos}
B.~Csik\'os, \emph{On the volume of the union of balls}, Discrete Comput. Geom. \textbf{20} (1998), 449--461, \href{https://doi.org/10.1007/PL00009395}{doi:10.1007/PL00009395}.

\bibitem{Gorbovickis}
I.~Gorbovickis, \emph{Kneser--Poulsen conjecture for a small number of intersections}, Contrib. Discrete Math. \textbf{9} (2014), no.~1, 1--10, \href{https://doi.org/10.55016/ojs/cdm.v9i1.62082}{doi:10.55016/ojs/cdm.v9i1.62082}.

\bibitem{GM}
Y.~Gordon and M.~Meyer, \emph{On the volume of unions and intersections of balls in Euclidean space}, Geometric Aspects of Functional Analysis, Oper. Theory Adv. Appl., vol.~77, Birkh\"auser, Basel, 1995, pp.~91--101, \href{https://doi.org/10.1007/978-3-0348-9090-8_9}{doi:10.1007/978-3-0348-9090-8\_9}.

\bibitem{Kneser}
M.~Kneser, \emph{Einige Bemerkungen \"uber das Minkowskische Fl\"achenmass}, Arch. Math. \textbf{6} (1955), 382--390.

\bibitem{Poulsen}
E.~T.~Poulsen, \emph{Problem 10}, Math. Scand. \textbf{2} (1954), 346.

\bibitem{WP}
P.~Wu and Q.~Peng, \emph{The Kneser--Poulsen volume inequalities for $d+4$ balls}, preprint, 2026, \href{https://hexagonmath.org/2611.01327v1}{Hexagon:2611.01327v1}, \href{https://doi.org/10.5281/zenodo.24623358}{doi:10.5281/zenodo.24623358}.
\end{thebibliography}
\end{document}